\documentclass[10pt]{article}
\usepackage[tbtags]{amsmath}
\usepackage{epsfig, amstext,amssymb,amsthm,latexsym}
\allowdisplaybreaks[4]

\usepackage{amssymb,color}

\definecolor{c20}{rgb}{0.,0.7,0.}
\definecolor{c30}{rgb}{0.,0.,1.}
\definecolor{c40}{rgb}{1,0.1,0.7}
\definecolor{c50}{rgb}{1,0,0}
\definecolor{c60}{rgb}{1,0.7,0.1}
\definecolor{c70}{rgb}{1,0.6,0.01}
\definecolor{c80}{rgb}{1,0.3,0.01}
\definecolor{c90}{rgb}{0.6,0.3,0.1}

\def\cga#1{\textcolor{c60}{#1}}
\def\cgb#1{\textcolor{c40}{#1}}
\def\cgc#1{\textcolor{c70}{#1}}
\def\cgd#1{\textcolor{c80}{#1}}

\def\cga#1{#1}
\def\cgb#1{#1}
\def\cgb#1{#1}
\def\cgc#1{#1}
\def\cgd#1{#1}
\def\ccr#1{#1}

\newtheorem{theo}{Theorem}[section]
\newtheorem{sat}[theo]{Proposition}
\newtheorem{de}[theo]{Definition}
\newtheorem{lem}[theo]{Lemma}
\newtheorem{example}[theo]{Example}
\newtheorem{korr}[theo]{Corollary}
\newtheorem{remark}[theo]{Remark}

\newcommand{\neprop}[1]{{ Proposition \ref{#1}}}

\newcommand{\nekorr}[1]{{ Corollary \ref{#1}}}

\newcommand{\kb}[1]{\boldsymbol{#1}}
\newcommand{\vk}[1]{\kb{#1}}

\def\FRE{\mbox{Fr\'{e}chet }}

\def\x{\vk{x}}

\def\kal#1{{\cal{ #1}}}

\def\fracl#1#2{\biggr( \frac{#1}{#2} \biggl) }

\newcommand{\abs}[1]{\lvert #1 \rvert}
\newcommand{\Abs}[1]{ \Bigl \lvert #1 \Bigr \rvert}

\newcommand{\E}[1]{\mbox{\rm$\vk{E}$}\{#1\}}

\newcommand{\pk}[1]{\mbox{\rm$\vk{P}$} \{#1\} }

\newcommand{\R}{\!I\!\!R}

\newcommand{\N}{\!I\!\!N}
\newcommand{\inr}{\in \R}

\newcommand{\inn}{\in \N}

\newcommand{\ldot}{,\ldots,}

\newcommand{\limit}[1]{\lim_{#1 \to   \infty}}
\newcommand{\as}{ \stackrel{a.s.}{\to}}
\newcommand{\todis}{\stackrel{d}{\to}}
\newcommand{\toprob}{ \stackrel{p}{\to}}

\newcommand{\equaldis}{\stackrel{d}{=}}

\newcommand{\BQN}{\begin{eqnarray}}
\newcommand{\EQN}{\end{eqnarray}}

\newcommand{\BQNY}{\begin{eqnarray*}}
\newcommand{\EQNY}{\end{eqnarray*}}

\newcommand{\BS}{\begin{sat}}
\newcommand{\ES}{\end{sat}}

\newcommand{\BT}{\begin{theo}}
\newcommand{\ET}{\end{theo}}
\newcommand{\BL}{\begin{lem}}
\newcommand{\EL}{\end{lem}}

\newcommand{\BK}{\begin{korr}}
\newcommand{\EK}{\end{korr}}

\newcommand{\BD}{\begin{de}}
\newcommand{\ED}{\end{de}}

\newcommand{\BIT}{\begin{itemize}}
\newcommand{\EIT}{\end{itemize}}
\newcommand{\BDI}{\begin{description}}
\newcommand{\EDI}{\end{description}}
\newcommand{\BEN}{\begin{enumerate}}
\newcommand{\EEN}{\end{enumerate}}
\newcommand{\BC}{\begin{cases}}
\newcommand{\EC}{\end{cases}}

\newcommand{\QED}{\hfill $\Box$}

\newcommand{\HR}{\notag \\}

\newcommand{\IF}{\infty}

\newcommand{\COM}[1]{}

\def\N{{\cal N}}
\def\R{{\cal R}}
\def\Z{\kal{Z}}

\def\ORD#1#2{ X_{#1- #2+1: #1} }
\def\ORDa#1#2{ #1_{N (t)- #2+1: N(t)} }
\def\ORDb#1#2{ #1_{N (t)- #2+1: N(t)} }

\def\bay{\bar{b}_1(t)}
\def\aay{\bar{a}_1(t)}
\def\bby{\bar{b}_2(t)}
\def\aby{\bar{a}_2(t)}
\def\bay{b_1(t)}
\def\aay{a_1(t)}
\def\bby{b_2(t)}
\def\aby{a_2(t)}

\def\stpX{S_1(p,t)}
\def\stpY{S_2(q,t)}
\def\kZ{\kal{Z}}
\def\kammap{K_p}
\def\kammaq{K_q}
\def\kammai{K_i}
\def\kammaa{K_1}
\def\kammam{K_m}

\begin{document}
\title{On the asymptotic distribution\\ of \cga{certain} bivariate reinsurance treaties}

\author{Enkelejd Hashorva\thanks{
Tel: +41 31 384 5528  Fax +41 31 384 4572\newline
Email: enkelejd.hashorva@stat.unibe.ch, \ enkelejd.hashorva@Allianz-Suisse.ch}\\
Allianz Suisse Insurance Company\\
Laupenstrasse 27,  3001 Bern, Switzerland\\
\\
 and\\
\\
Institute of Mathematical  Statistics and Actuarial Science \\
University of Bern\\
Sidlerstrasse 5, 3012 Bern, Switzerland}

\bigskip

 \maketitle

{\bf Abstract:}
Let $\{(X_n,Y_n), n\ge 1\}$ be  bivariate random claim sizes with
common distribution function $F$ and let $\{N(t), t \ge 0\}$ be a stochastic
process which counts the number of claims that occur in the time interval $[0,t], t\ge 0$.
In this paper we derive the joint asymptotic  distribution of randomly indexed order statistics of the
random sample $(X_1,Y_1),(X_2,Y_2)\ldot (X_{N(t)},Y_{N(t)})$ which is then used to obtain
asymptotic representations for the joint distribution of  two generalised largest claims reinsurance
treaties available under specific insurance settings. \ccr{As a by-product we obtain
a stochastic representation of a $m$-dimensional $\Lambda$-extremal variate in terms of iid unit
exponential random variables.}

\noindent{\em \bf Key Words:} Generalised largest claims reinsurance treaty,
bivariate random order statistics, \cga{bivariate treaties}, asymptotic results, extreme value theory,
$m$-dimensional $\Lambda$-extremal variate.

{\it MSC index}: 60G70 (Primary) 
     62P05 (Secondary)


\section{Introduction}
Let $\{(X_n,Y_n), n\ge 1\}$ be  bivariate random claim sizes with
common distribution function $F$ arising from an insurance
portfolio. In a specific insurance
context, $X_n$ may stay for instance for the total claim amount related to  the
$n$th accident, and $Y_n$ for the corresponding total expense amount.
Let  further $\{N(t), t \ge 0\}$ be a stochastic process which counts the number of claims that
occur in the time interval $[0,t], t\ge 0$. So we observe  $(X_1,Y_1)  \ldot
(X_{N(t)}, Y_{N(t)})$ claims up to time $t>0$. Denote by $X_{i: N(t)},Y_{i: N(t)}, 1 \le i \le N(t)$ the
corresponding  $i$th \cgb{lowest} claim in the above random sample \cgb{taken component-wise}.\\
Especially in reinsurance applications, the largest claims
$\ORDa{X}{i},\ORDb{Y}{i}, i=1 \ldot m$ are of particular importance.
For instance, consider the following reinsurance contract introduced
by Ammeter (1964)
\BQNY  \stpX &=& X_{N(t):N(t)}+ X_{N(t)-1:N(t)}+
\cdots + X_{N(t)-p+1:N(t)},\quad p\ge 1
\EQNY
which is known as the Largest Claim Reinsurance. Thus at the time point $t$ a reinsurer
covers the total loss  amount $\stpX$ which is the simplest linear transformation of the upper
$p$ largest claims observed up to time $t$.  Another simple reinsurance treaty
\BQNY    \stpX &=& X_{N(t):N(t)}+ X_{N(t)-1:N(t)}+ \cdots + X_{N(t)-p+2:N(t)} - (p-1)X_{N(t)-p+1:N(t)}, \quad p\ge 2
\EQNY
is introduced by Th{\'e}paut (1950). This treaty is known in actuarial literature as the ECOMOR
reinsurance treaty (for further details  \cgc{see e.g.\
Teugels (1985), Beirlant et al.\ (1996), Embrechts et al.\ (1997), Rolski et al.\ (1999), and Mikosch (2004)}).

In general one can define the total loss amount by
\BQNY   \stpX &=&\sum_{j=1}^p g_j^*(\ORDa{X}{j}),
\EQNY
with $g_j^*, 1 \le j\le p$  real measurable functions. Typically, the prime
interest of a reinsurer is the calculation of the pure premium
$\E{\stpX}$. This can be clearly expressed as  $\sum_{j=1}^p
\E{g_j^*(\ORDa{X}{j})}$ supposing  additionally that
$g_j^*(\ORDa{X}{j}), 1 \le j\le p$ have finite expectations.
\cgb{It is often assumed in actuarial applications} that $N(t)$ is
independent of claim sizes. With that assumption we have for each
$j=1 \ldot p$
\BQNY\E{ g_j^*(\ORDa{X}{j})} &=&  \E{ \E{ g_j^*(\ORDa{X}{j})\lvert N(t)}}
=\sum_{k=1}^\IF  \E{g_j^*(\ORD{k}{j})}\pk{N(t)=k}.
\EQNY
Essentially, in order to compute the pure premium, we need to calculate the expectations of  the upper order
statistics. Further we need to know the distribution function of the counting random variable $N(t)$.
Next, since we consider a bivariate setup, let us suppose further that
a similar treaty
\BQNY   \stpY &=&\sum_{j=1}^q g_j^{**}(\ORDa{Y}{j}), \quad q\ge 1,
\EQNY
with $g_j^{**}, 1 \le j\le q$ real measurable functions, covers the risks
modelled by $Y_1,Y_2 \ldot Y_{N(t)}$.

In order to price both \cgb{treaties mentioned above},
the reinsurer needs to have some indications concerning the distribution function
of the total loss amount $(\stpX,\stpY)$. For instance, if in particular  the standard deviation (variance)
premium principle is used, then an estimate of the standard deviation is required.
Further the dependence \cgc{between} $\stpX$ and $\stpY$ needs to be quantified.

An asymptotic model for both treaties can be regarded as a good candidate
to overcome the difficulties in specification of
the model (\cgc{distribution assumptions for the claim sizes  or assumptions
on the first and second moment of $X_1,Y_1$)}. The idea is
to let $t\to \infty$ and to investigate the joint asymptotic behaviour of  $(\stpX, \stpY)$.

Relying on extreme value theory  Embrechts et al.\ (1997) (see Example 8.7.7 therein)
derives the asymptotic distribution of both the Largest Claim Reinsurance and the ECOMOR Reinsurance treaties
when $N(t)/t \to \lambda \in (0,\infty) $ in probability. In Hashorva (2004) the asymptotic limiting
distribution for $(\stpX , \stpY)$ a bivariate ECOMOR reinsurance treaty is obtained when the marginal distributions
of $F$ are in the max-domain of attraction of \cgd{the} Gumbel distribution, imposing further the iid
(independent and identically distributed) assumption on the claim sizes.

With main impetus from the afore-mentioned results, we consider in this article
special reinsurance treaties with $g_j^*, g_j^{**}$
simple linear functions and claim sizes which can be dependent (thus dropping the iid assumption).\\
In Section 2 we  first deal with the joint distribution function of randomly indexed upper
order statistics; application to reinsurance and details for the iid case are presented in the Section 3.
Proofs \ccr{of the results} are given in Section 4.

\section{Joint limiting distribution of randomly indexed\\
 upper order statistics}

It is well-known (see e.g.\ Reiss (1989)) that the joint
convergence in distribution under the assumption of iid claim
sizes
\BQN \label{eq:1} \lefteqn{\Biggl(
\biggl(\frac{X_{n:n}-b_1(n)}{a_1(n)},\frac{Y_{n:n}-b_2(n)}{a_2(n)}\biggr)\ldot
\biggl(\frac{X_{n-m+1:n}-b_1(n)}{a_1(n)},\frac{Y_{n-m+1:n}-b_2(n)}{a_2(n)}\biggr)\Biggr) }\hspace{14 cm} \notag \\
\todis \Bigl( (\kal{X}_1, \kal{Y}_1)\ldot (\kal{X}_m,
\kal{Y}_m)\Bigr), \quad n\to \IF \hspace{5 cm}
\EQN
with $(\kal{X}_i,\kal{Y}_i),i\le m$ random vectors holds
for all integers $m\inn$ and given real functions $a_i(t)>0, b_i(t),  i=1,2$ iff the
underlying distribution function $F$ of the random vector
$(X_1,Y_1)$ is in the max-domain of attraction of  $H$, the
 max-stable bivariate distribution function of $(\kal{X}_1, \kal{Y}_1)$, i.e.\
\begin{align}\label{eq:weak2}
 \limit{t} \sup_{(x,y)\inr^2}\Abs{F^t(a_1(t)x+\cgb{b_1}(t), a_2(t)y+b_2(t)) - H(x,y)} &=& 0
\end{align}
holds. See below \eqref{eq:hm} for the joint distribution of $\kal{X}_1\ldot \kal{X}_k,k\ge 1$.
 For short we write the above fact as $F\in MDA(H)$. \\
We note in passing that the standard notation $\todis, \toprob ,
\as $ which we use throughout in this paper mean convergence in
distribution,   convergence in probability  and almost sure
convergence, respectively.

Actually \eqref{eq:weak2} implies that for the marginal distributions of $F$ (denoted here by $F_1,F_2$) we have $F_i\in MDA(H_i),i=1,2$, where
the standard extreme value distribution function $H_i$ is either of the following
\def\al{\alpha}
\BQNY \cga{\Phi_{\al_i}}(x) &=&\exp(-x^{-\alpha_i}), \quad \alpha_i>0, \quad x>0,
\EQNY
or
\BQNY   \Psi_{\al_i}(x) &=&\exp(-\abs{x}^{\alpha_i}), \quad \alpha_i>0, \quad x<0,
\EQNY
or
\BQNY   \Lambda(x) &=& \exp(-\exp(-x)), \quad x \inr,
\EQNY
i.e.\   unit \cgb{Fr\'{e}chet}, Weibull or Gumbel  distribution, respectively.

\ccr{For iid claim sizes}, if $N(t) \toprob \infty$ as $t\to \infty$ and $N(t)$ is
independent of the claim sizes for all $t$ large, then by Lemma \cgd{2.5.6} of Embrechts et al.\ (1997) the asymptotic
relation \eqref{eq:weak2} implies for any $i\ge 1$
\BQNY
\Biggl(\frac{X_{N(t)-i+1:N(t)}-b_1(N(t))}{a_1(N(t))},\frac{Y_{N(t)-i+1:N(t)}-b_2(N(t))}{a_2(N(t))}\Biggr)
&\todis& (\kal{X}_i, \kal{Y}_i), \quad t\to \IF.
\EQNY
A different situation arises when transforming with $a_i(t),b_i(t)$ instead of
the random functions \\
$a_i(N(t)),b_i(N(t)), i=1,2$. So, if we assume further that
\BQN \label{eq:N} \frac{N(t)}{t} & \toprob & \kal{Z}, \quad t \to \IF
\EQN
holds with $\kal{Z}$ such that $\pk{\kal{Z}>0}=1$, then it follows (along the lines of Theorem
4.3.4 of Embrechts et al.\ (1997)) that for all $i\ge 1$
\BQNY
\biggl(\frac{X_{N(t)-i+1:N(t)}-b_1(t)}{a_1(t)},\frac{Y_{N(t)-i+1:N(t)}-b_2(t)}{a_2(t)}\biggr)
&\todis & (\kal{X}^*_i, \kal{Y}^*_i), \quad t\to \infty
\EQNY
holds with $(\kal{X}^*_i, \kal{Y}^*_i)$ a new bivariate random vector.\\
\ccr{For the univariate case Theorem 4.3.4 of Embrechts et al.\ (1997)  shows
the distribution functions of $\kal{X}^*_i$.
Proposition 2.2 of Hashorva (2003) gives an explicit expression
for the joint distribution of $(\kal{X}^*_1 \ldot \kal{X}^*_i).$
}
\ccr{Our next result is more general. We consider the bivariate setup allowing claim sizes to be dependent, and further instead of
\eqref{eq:N} we assume the convergence in distribution
\BQN \label{eq:NN} \frac{N(t)}{t} &\todis &\kal{Z}, \quad t \to \IF,
\EQN
}
with $\kal{Z}$ almost surely positive. Clearly the above condition is weaker than \eqref{eq:N}.

\BS\label{prop:1}
Let  $\{(X_n,Y_n), n\ge 1\}$ be  bivariate claim sizes  with common distribution function $F$
independent of the counting process $\{N(t), t \ge 0\}$ for all $t\ge 0$.
If condition \eqref{eq:NN} is fulfilled with $\kal{Z}$ positive and non-zero and
 further  \eqref{eq:1} holds for some fixed $m\inn$ with constants $a_i(t)>0,b_i(t),i=1,2$
that satisfy \eqref{eq:weak2}, then we have for  $t \to \IF$
\BQN\label{proof1}
\lefteqn{\Biggl( \biggl(\frac{X_{N(t):N(t)}- \bay }{ \aay} ,\frac{Y_{N(t):N(t)}-  \bby }
{ \aby}\biggr)\ldot
\biggl(\frac{X_{N(t)-m+1:N(t)}- \bay }{\aay} ,\frac{Y_{N(t)-m+1:N(t)}- \bby }{\aby}\biggr) \Biggr)}\notag\\
&\todis& \Bigl( (\kZ^{\gamma_{1}}\kal{X}_1+ \delta_{1} ,\kZ^{\gamma_{2}}\kal{Y}_1+ \delta_{2})
\ldot (\kZ^{\gamma_{1}}\kal{X}_m+\delta_{1}, \kZ^{\gamma_{2}}\kal{Y}_m+\delta_{2})\Bigr),
\hspace{5cm}
\EQN
with \cgd{$ \gamma_{i}:=1/\alpha_{i}, -1/ \alpha_{i}, 0$} if
$\cgd{F_i \in MDA(\Phi_{\alpha_i}), MDA(\Psi_{\alpha_i}), \ } MDA(\Lambda)$, respectively.
Further $\delta_{i}:= \ln \kal{Z}$ if $F_i \in MDA(\Lambda)$ and 0 otherwise for i=1,2.
\ES
In the above proposition we do not assume explicitly the independence of the claim sizes.
 For iid claim sizes 
(recall that \eqref{eq:weak2} is equivalent to \eqref{eq:1}) we obtain  immediately:

\BK\label{korr:1} Let  $\{(X_n,Y_n), n\ge 1\}$ be  iid
bivariate claim sizes  with distribution function $F$ independent of $N(t), t\ge 0$.
If condition \eqref{eq:weak2} is satisfied and further \eqref{eq:NN} holds with $\kal{Z}$ almost surely positive,
then  \eqref{proof1} holds for any $m\ge 1$ with  $\delta_i, \gamma_i,i=1,2$ as in \neprop{prop:1}.
\EK

 {\bf Remarks:} (i) If $\{X_n, n\ge 1\}$ are iid with
distribution function $F\in MDA(H)$, with $H$ a univariate extreme
value distribution, then the joint density function of
$(\kal{X}_1 \ldot \kal{X}_m), m\ge 1$ is given by \ccr{(see e.g.\ Embrechts et al.\ (1997) p.\ 201)}
\BQN\label{eq:hm}
h_m(\x)&=&H(x_m) \prod_{i=1}^m \frac{H'(x_i)}{H(x_i)}, \quad
\text{with } \cgb{x_1 > x_2 >\cdots > x_m},\quad  \prod_{i=1}^m
H(x_i)\in (0,1).
\EQN
\ccr{Referring to Embrechts et al.\ (1997) the random vector $(\kal{X}_1 \ldot \kal{X}_m)$ is called
a $m$-dimensional $H$-extremal variate. If} $H$ is the unit Gumbel distribution,
then the following stochastic representation
(see Theorem 7.1 of Pakes and Steutel (1997))
\BQN\label{eq:stoch:rep}
\Bigl(\kal{X}_i-\kal{X}_{i+1}\Bigr)_{i=1 \ldot k}  \equaldis
\Bigl( \frac{E_i}{i}\Bigr) _{i=1 \ldot k}, \quad k\ge 1
\EQN
holds with $E_i, i\ge 1$ iid unit exponential random variables. The standard notation $\equaldis$ means
equality of distribution functions. \ccr{See \eqref{eq:stoch:L} for
a stochastic representation of $(\kal{X}_1 \ldot \kal{X}_m)$.}

(ii) In the above proposition $\kal{Z}$ is independent of $(\kal{X}_i,\kal{Y}_i), 1 \le i \le m$.
This follows immediately recalling that the claim sizes are independent of the counting process $N(t),t\ge 0$.

For the more general class of strictly stationary random sequences
conditions for joint weak convergence of upper order statistics are available in literature.
For  $\alpha$-mixing stationary sequences (univariate case)  the mentioned  weak convergence
is discussed in Hsing (1988). In Theorem 4.1 and Theorem 4.2 therein necessarily and sufficient
conditions for joint weak convergence are derived. Moreover, it is shown that the limit distribution should
have a specific form. Novak \cgb{(2002)} gives a simplified version of the limit distribution.
A variant of Theorem 4.2 of Hsing (1988) can be found in Leadbetter (1995).
Convergence results for dependent bivariate random sequences are obtained in H\"usler (1990,1993).

\section{Joint Asymptotic Distribution of $\stpX$ and  $\stpY$}
In this section we investigate the asymptotic behaviour ($t \to \infty$) of $(\stpX, \stpY)$ \ccr{defined in the introduction.}
In general some restrictions on the choice  of the functions $g_j^*, 1 \le j\le p,g_j^{**}, 1 \le j\le q$
should be imposed. Obviously the total expected  loss should be non-negative.
Tractable (simple) functions which we consider here are
\BQN \label{eq:gy}
 g_j^*(x) = k_{j1}x,\quad 1 \le j \le p, \quad \text{ and }   g_j^{**}(x) =k_{j2}x ,\quad 1 \le j\le q,
\EQN
with $k_{j1},1 \le j\le p, k_{j2}, 1 \le j\le q$ real constants.  Put throughout in the following
$c_1:= \sum_{j=1}^p k_{j1}, c_2:= \sum_{j=1}^q k_{j2}$. Based on the previous results,
we consider now the joint asymptotic behaviour of the reinsurance treaties.
\BS\label{prop:2}
Let $\{(X_n,Y_n), n\ge 1\}$,   $\{N(t),  t \ge 0\}, \kal{Z}$ and $\gamma_{i}, \delta_{i}, i=1,2$ be as
in \neprop{prop:1}, and for $p,q\inn$ let $g_j^*, 1 \le j\le p,g_j^{**}, 1 \le j\le q$  be as in \eqref{eq:gy}.
If condition \eqref{eq:N} is fulfilled and \eqref{eq:1} holds for some $m\ge \max(p,q)$, then
\BQN \label{eq:JHa}
\Biggl(\frac{\stpX-  \bay c_1  }{\aay},\frac{\stpY-  \bby c_2 }{\aby}\Biggr)
& \todis &
\Biggl(  \kZ^{\gamma_{1}} \sum_{j=1}^{p}  k_{j1} \kal{X}_j+ c_1\delta_{1} ,
\kZ^{\gamma_{2}} \sum_{j=1}^{q}  k_{j2}\kal{Y}_j + c_2\delta_{2} \Biggr)
\EQN
holds as $t\to \IF$.
\ES

{\bf Remarks:} (i) In the above proposition there is no
restriction imposed on the constants $k_{ji}.$ Neither need the
claim sizes be positive. Typically, in reinsurance claim sizes are assumed to be positive, and further, some restrictions
have to be imposed on the constants $k_{ji}$. Implicitly we may
require that these constants are such that the expected loss (at
any time point) for both $S_1(p,t), S_2(q,t)$ is non-negative. A
more explicit assumption would be to suppose that both $c_1,c_2$ are positive.

(ii) \neprop{prop:1} and \neprop{prop:2} can be shown with similar arguments for
a general multivariate setup, i.e.\ for claim sizes being random vectors in $\R^k, k\ge 3$.

We discuss next the  bivariate counterpart of Example 8.7.7
of Embrechts et al.\ (1997), which was our starting point.
We correct a missing constant in an asymptotic result therein.

In order to keep things simple, we impose throughout in the following  the assumptions of
\nekorr{korr:1}, considering therefore iid claims sizes
independent of the counting process $N(t)$.

Note in passing that in Example 8.7.7 of Embrechts et al.\ (1997)
the claim sizes (univariate setup) are iid being further
independent of the counting process, which satisfies \eqref{eq:N}
with $\kal{Z}=\lambda \in (0,\infty)$ almost surely.

To this end, suppose for simplicity that both marginal distributions $F_1,F_2$ of $F$ are identical
and let $F_1^{-}, F_2^{-}$ be the generalised inverse of $F_1$ and $F_2$, respectively.

{\bf Case a)}: If $F_1\in MDA(\Phi_{\alpha}), \alpha> 0$, then we may take  $b_1(t)=b_2(t)=0$ and
$$a_1(t)=a_2(t):=\inf\{x\inr: F_1(x)> 1- 1/t\}= F^{-}_1(1- 1/t), \quad \ccr{t>1}$$
(see e.g.\ Resnick (1987), Reiss (1989) or Embrechts et al.\ (1997)).
So we obtain for $p, q \inn$ and \cgd{$k_{ji}$ real constants}
\BQN\label{eq:C:1}
\Biggl(\frac{\stpX}{a_1(t)},
\frac{\stpY  }{a_2(t)} \Biggr)
& \todis & \Biggl(  \kal{Z}^{1/\alpha} \sum_{j=1}^{p}  k_{j1} \kal{X}_j ,
\kal{Z}^{1/\alpha} \sum_{j=1}^{q}  k_{j2}\kal{Y}_j \Biggr)
\EQN
as $t\to \infty$.
\cgd{Embrechts et al.\ (1997) proves for the ECOMOR case
$$\frac{ S_1(p,n)}{ a_1(n)} \todis   \lambda ^{1/\alpha}
\sum_{i=1}^{p-1} i  (\kal{X}_i- \kal{X}_{i+1} ),
\quad
n\to \infty $$
which follows immediately from \eqref{eq:C:1} putting $\kal{Z}=\lambda>0$.
}

{\bf Case b)}:
When $F_1\in MDA(\Psi_{\alpha}), \alpha > 0$ holds, then we may take
$b_1(t)=b_2(t):=\omega, t>1,$ with $\omega:=\sup\{x: F_1(x)<1\}$ the upper endpoint of the distribution
function $F_1$ which is necessarily finite and define
$$a_1(t)=a_2(t):=\omega - F^{-}_1(1- 1/t),\quad  t>1.$$
Thus we have
\BQNY
\Biggl(\frac{\stpX -c_1\omega }{a_1(t)},
\frac{\stpY -c_2\omega  }{a_2(t)}\Biggr)
& \todis & \Biggl( \kal{Z}^{-1/\alpha} \sum_{j=1}^{p}  k_{j1} \kal{X}_j ,
\kal{Z}^{-1/\alpha} \sum_{j=1}^{q}  k_{j2}\kal{Y}_j \Biggr),  \quad  t \to \infty.
\EQNY

{\bf Case c)}: If $F_1\in MDA(\Lambda)$, then we put for $t$ large
$$ b_1(t)=b_2(t):= F^{-}_1(1- 1/t)$$
and
$$ a_1(t)=a_2(t):=\int_{b_1(t)}^{\omega} [1- F_1(s)]/[1- F_1(b_1(t))]\, d s.$$
Thus we have that the \cga{right hand side} of \cgb{\eqref{eq:JHa}} is given by
\BQNY
\Biggl( \sum_{j=1}^{p}  k_{j1} \kal{X}_j +c_1\ln \kal{Z},
\sum_{j=1}^{q}  k_{j2}\kal{Y}_j+c_2\ln \kal{Z}\Biggr).
\EQNY

Next, we consider 3  examples.

{\bf Example 1.} (Generalised ECOMOR reinsurance treaty).  Assume that the constants $k_{ij}$
are such that $c_1=c_2=0.$ This is fulfilled in the special case
of the ECOMOR treaty (see introduction above).

Thus under the assumptions of Proposition \ref{prop:2} we have for $F_1,F_2\in MDA(\Lambda)$ the
convergence in distribution
\BQNY
\Biggl(\frac{\stpX }{\aay},
\frac{\cgc{\stpY }}{\aby}\Biggr)
& \todis & \Biggl( \sum_{j=1}^{p}  k_{j1} \kal{X}_j ,
\sum_{j=1}^{q}  k_{j2}\kal{Y}_j  \Biggr), \quad t\to \infty.
\EQNY
\cgd{For the ECOMOR treaty Embrechts et al.\ (1997) obtains \ccr{with $\kal{Z}=\lambda$ almost surely} 
\BQNY
\frac{ S_1(p,n)}{ a_1(n)} &\todis &   \sum_{j=1}^{p-1} j (\kal{X}_j- \kal{X}_{j+1})
 \equaldis \sum_{j=1}^{p-1}E_j, \quad n\to \infty,
\EQNY
with $E_j, 1 \le j\le p-1$ iid unit exponential random variables, which follows
also by the bivariate result above (recall \eqref{eq:stoch:rep}).
}

There is a remarkable fact  in the above asymptotic result, namely the random variable $\Z$ does not appear in the right hand side of the asymptotic expression.
This is not the case in general for $F_i,i=1,2$ in the max-domain of attraction of \FRE or Weibull.

{\bf Example 2.} (Asymptotic independence of the components of the maximum claim sizes.)
As mentioned in the introduction, for the asymptotic
considerations, we are not directly interested in the joint
distribution function $F$ of the claim sizes, but on the limiting distribution $H$.
In some applications, even if $F$  is not a product distribution, it may happen that  $H$ is a product
distribution, meaning $H(x,y)=H_1(x)H_2(y), \forall x,y\inr.$ This
implies that $\kal{X}_i$ is independent of $\kal{Y}_j$ for any
$i,j\ge 1$ and thus $[\stpX -b_1(t)c_1 ]/\aay$ is asymptotically
independent of $[\stpY -b_2(t)c_2]/\aby$ as $t\to \infty$. \cgd{So the asymptotic distribution of each treaty can be easily calculated using \eqref{eq:hm}.
We note in passing that there are several known conditions for asymptotic independence, see e.g.\
 Galambos  (1987), Resnick (1987), Reiss (1989),  Falk et al.\ (1994),
H\"usler (1994). 
}

{\bf Example 3.} ($F_1,F_2$ with exponential tails and $N(t)$ Poisson).
Consider iid claim sizes with joint distribution function $F$ which has marginal distributions
tail equivalent to the unit exponential distribution function, i.e.
$ \lim_{x \to \infty} \exp(x)[1-F_i(x)]=1, i=1,2$.
\cgc{It follows that $F_1,F_2 \in MDA(\Lambda)$} with
constants $a_1(t)=a_2(t)=1, b_1(t)=b_2(t)=\ln t, t>0$.
As in the above example, assume further that the distribution function $H$
is a product distribution and \cgb{$N(t)$ is a (homogeneous) Poisson process with
parameter $\lambda>0$ independent of the claim sizes.}
We have thus $N(t)/ t \as  \lambda$
and for any $p,q\inn$
\BQNY
\Biggl(\stpX - c_1 \ln t, \stpY - c_2 \ln t \Biggr)
& \todis & \Biggl( \sum_{j=1}^p k_{j1} \kal{X}_{j} + c_1 \ln \lambda ,
\sum_{j=1} ^q k_{j2} \kal{Y}_{j}+ c_2 \ln \lambda \Biggr),  \quad t\to \infty,
\EQNY
with $\kal{X}_j$ independent of $\kal{Y}_j$.  \ccr{Recall $c_1:=\sum_{j=1}^p k_{j1}, c_2:=\sum_{j=1}^q k_{j2}$.} \\
Borrowing the idea of Example 8.7.7 of Embrechts et al.\ (1997) we  find an explicit  formula for the right hand side above.
Let therefore $E_{ji}, j\ge 1, i=1,2$ be iid unit exponential random variables and put
 $ \bar k_{li}:=\sum_{j=1}^l k_{ji}/l, l\ge 1, i=1,2$. It is well-known that
(see e.g.\ Reiss (1989))
$$\Bigl(E_{n-j+1:n,i}\Bigr)_{j=1 \ldot n} \equaldis
  \Bigl(\sum_{l=j}^n \frac{E_{li}}{l}\Bigr)_{j=1 \ldot n}, \quad n\ge 1, i=1,2.
 $$
\ccr{It is well-known that for large $n$
$$ \sum_{l=1}^n \frac{1}{l} - \ln n= K+o(1)$$
where $K$ is the Euler-Mascheroni constant.} So we may write for \ccr{$n\ge p\ge 1$}
\BQNY
\sum_{j=1}^p k_{j1} ( E_{n-j+1:n,1}- \ln n ) &\equaldis &
\sum_{j=1}^p k_{j1} \sum_{l=j}^n \frac{E_{l1}}{l}- c_1\ln n \\
&=& \sum_{l=1}^p  E_{l1} \frac{1}{l} \Bigl(\sum_{j=1}^l k_{j1} \Bigr)+
  c_1 \sum_{l=p+1}^n \frac{E_{l1}}{l}+c_1\Bigl( K- \sum_{l=1}^{n} \frac{1}{l} \Bigr)+o(1)\\
&\todis &  \sum_{l=1}^p  \bar k_{l1} E_{l1}  + c_1 \sum_{l=p+1}^\infty
\frac{E_{l1}-1}{l}+c_1 \Bigl( K- \sum_{l=1}^{p} \frac{1}{l} \Bigr),  \quad n\to \infty.
\EQNY
On the other hand
\BQNY
(E_{n:n,1}- \ln n \ldot  E_{n-j+1:n,1}- \ln n)  &\todis & (\kal{X}_{1}\ldot \kal{X}_{j}), \quad n\to \infty,
\EQNY
hence the continuous mapping theorem (see e.g.\ Kallenberg (1997)) implies
\BQNY
 \sum_{j=1}^p k_{j1} ( E_{n-j+1:n,1}- \ln n ) & \todis& \sum_{j=1}^p k_{j1} \kal{X}_{j}, \quad n\to \infty.
 \EQNY
Thus we have the stochastic representation
\BQN\label{eq:stochA}
\sum_{j=1}^p k_{j1} \kal{X}_{j} &\equaldis  & \sum_{l=1}^p  \bar k_{l1} E_{l1}  + c_1 \sum_{l=p+1}^\infty
\frac{E_{l1}-1}{l}+c_1\kammap,
\EQN
with $\kammai:=K- \sum_{l=1}^{i} \frac{1}{l}, i\ge 1.$ We obtain  thus proceeding similarly for
the second treaty and recalling the asymptotic independence assumption
\BQNY
\lefteqn{\Biggl(\stpX - c_1 \ln t, \stpY - c_2 \ln t \Biggr)}\\
& \todis &
\Biggl( \sum_{j=1}^{p}  \bar{k}_{j1} E_{j1}+ c_1 \Bigl[
\sum_{j=p+1}^{\infty} \frac{E_{j1}- 1}{j}+ \kammap +\ln \lambda\Bigr]
  ,\sum_{j=1}^{q}  \bar{k}_{j2} E_{j2}+ c_2 \Bigl[
  \sum_{j=q+1}^{\infty} \frac{E_{j2}- 1}{j}+\kammaq+\ln \lambda\Bigr] \Biggr).
\EQNY

{\bf Remarks.}   (i) For a suitable choice of constants
\eqref{eq:stochA} implies for any $m\ge 2$ \BQN\label{eq:stoch:L}
(\kal{X}_{1} \ldot \kal{X}_{m}) \equaldis   \Bigl( E_{11}  +
\sum_{l=2}^\infty \frac{E_{l1}-1}{l}+\kammaa\ldot E_{m1}  +
\sum_{l=m+1}^\infty \frac{E_{l1}-1}{l}+\kammam\Bigr), \EQN hence in
particular we have for any $i\in \N$
$$\E{\kal{X}_{i} }= 1+\kammai, \quad {\bf Var}\{\kal{X}_{i} \}=
1+ \sum_{l=i+1}^\infty \frac{1}{l^2} = \frac{\pi^2}{6}+1- \sum_{l=1}^i \frac{1}{l^2}<\infty . $$
(ii) In Example 8.7.7 of Embrechts et al.\ (1997) the term  $-k \ln k$ in page 519 should be
$k(K- \sum_{j=1}^{k} \frac{1}{j}).$

A general result is given in the following proposition.

\BS\label{prop:3} Let $\{(X_n,Y_n), n\ge 1\}$ be iid bivariate
claim sizes with distribution function $F$. Assume that \eqref{eq:weak2} holds with given functions $a_i(t)>0, b_i(t), i=1,2$ and
$H(x,y)= \exp(-\exp(-x)-\exp(-y)), x,y\inr$. Suppose further that $\{N(t),  t \ge 0\}$  is a counting process independent
of the claim sizes satisfying \eqref{eq:NN} with $\kal{Z}$ positive and non-zero. Then we have the convergence in distribution
\BQN
\lefteqn{\Biggl(\frac{\stpX-  \bay c_1  }{\aay},
\frac{\stpY-  \bby c_2 }{\aby}\Biggr)}\notag \\
& \todis &
\ccr{\Biggl(\sum_{j=1}^{p}  \bar{k}_{j1} E_{j1}+ c_1 \Bigl[
\sum_{j=p+1}^{\infty} \frac{E_{j1}- 1}{j}+\ln \kal{Z} + \kammap\Bigr]
  ,\sum_{j=1}^{q}  \bar{k}_{j2} E_{j2}+ c_2 \Bigl[ \sum_{j=q+1}^{\infty} \frac{E_{j2}- 1}{j}+\ln \kal{Z} + \kammaq\Bigr] \Biggr)
 \notag
 }\\
\EQN as $t\to \infty$, with  $c_i,  \bar k_{ji}, \kammap,\kammaq$  as
defined above and $E_{ji}, j\ge 1, i=1,2$ iid unit exponential random variables.
 \ES

\section{Proofs}

{\bf Proof of Proposition \ref{prop:1}:} Define in the following for $t$ positive
$$
\ccr{\vk{T}_j(N(t))}:=(T_{j1}(N(t)),T_{j2}(N(t)))
=
\Biggr(
\frac{X_{N(t)-j+1:N(t)}-b_1(N(t))}{a_1(N(t))},
\frac{Y_{N(t)-j+1:N(t)}-b_2(N(t))}{a_2(N(t))}\Biggr)
$$
and 
$$
Z_t:=\frac{N(t)}{t}, \quad
\tilde{a}_i(t,z):= \frac{a_i( t z )}{a_i(t)},  \quad
\tilde{b}_i(t,z):=\frac{b_i(tz)-b_i(t) }{a_i(t)}, \quad z>0,
i=1,2. $$ Let $z_t, t\ge 0$ be arbitrary positive constants such
that $\limit{t} z_t= z\in (0,\infty)$. Assumption \eqref{eq:weak2} implies that the positive norming
functions $a_1,a_2$ are regularly varying (see e.g.\ Resnick (1987)),
hence using Theorem A3.2 of Embrechts et al.\ (1997) we obtain
$$ \limit{t} \tilde{a}_i( t,z_t)=z^{\gamma_i}, \quad i=1,2, $$
with $\gamma_i:= 1/ \alpha_i, -1/ \alpha_i, 0$
if $F_i\in MDA(\Phi_{\alpha_i}), F_i\in MDA(\Psi_{\alpha_i})$ or
$F_i\in MDA(\Lambda)$, respectively.\\
If $F\in MDA(\Psi_\alpha)$ or
$F\in MDA(\Phi_\alpha)$ with $\alpha>0,$ then for large $t$ we have $b_i(ts)-b_i(t)=0,i=1,2, \forall s>0$.
It is well-known (see Proposition 0.10, Proposition 1.1, Corollary 1.7 and Exercise 0.4.3.2 of Resnick (1987)) that
if  $F\in MDA(\Lambda)$ then  $\limit{t}(b_i(ts)-b_i(t))/a_i(t)= \ln s,  s>0$ holds locally unifromly
on $(0,\infty)$. Consequently we have
$$\limit{t} \tilde{b}_i(t,z_t)= \delta_i\ln z ,\quad i=1,2,$$
with $\delta_i=1$ if  $F_i\in MDA(\Lambda)$ and $\delta_i=0$, otherwise.

\ccr{In view of \eqref{eq:NN} applying now Theorem 3.27 of Kallenberg (1997) we obtain  the convergence in distributions
$$ \tilde{a}_i(t, Z_t) \todis \kal{Z}^{\gamma_i}, \quad \text{and  }
\tilde{b}_i(t,Z_t)\todis \delta_i\ln \kal{Z} ,\quad t\to \infty, i=1,2,
$$
}
hence as  $t \to \infty$
\BQNY
 \vk{U}_{t,N(t)} :=( \tilde{a}_1(t,Z_t),\tilde{a}_2(t,Z_t) ,\tilde{b}_1(t,Z_t),\tilde{b}_2(t,Z_t) )
&\ccr{\todis} & (       \kal{Z}^{\gamma_1},\kal{Z}^{\gamma_2},  \delta_1 \ln \kal{Z},  \delta_2 \ln \kal{Z})=:\vk{U}_{\kal{Z}}.
\EQNY
\ccr{Let $\vk{w}_j \inr^2, 1 \le j\le m$ and $\vk{u} \in (0,\infty)^2\times \R^2$ be given constants.
Since the claim sizes are independent of $N(t), t\ge 0$ we obtain by conditioning}
\ccr{
\BQNY \label{eq:uf}
\pk{\vk{U}_{t,N(t)} \le \vk{u}, \vk{T}_{1}(N(t))\le \vk{w}_1 \ldot \vk{T}_{m}(N(t))
 \le \vk{w}_m  \lvert N(t) =n } &=&\vk{1}( \vk{U}_{t,n} \le \vk{u}) h(n), \quad n\ge m, t>0,
\EQNY
with $h(n):=\pk{ \vk{T}_{1}(n) \le  \vk{w}_1 \ldot \vk{T}_{m}(n)\le \vk{w}_m}$.}
Here  $\vk{a} \le \vk{b}, \vk{a},\vk{b} \inr^k,k \ge 2$ is understood component-wise  and $\vk{1}(\cdot)$ is the indicator function.
\ccr{By the assumptions
\BQNY
 \limit{n} h(n)
&=&\pk{(\kal{X}_1, \kal{Y}_1)  \le  \vk{w}_1 \ldot   (\kal{X}_m, \kal{Y}_m ) \le \vk{w}_m}.
\EQNY
Using \eqref{eq:NN} we obtain  $N(t) \toprob \infty$ as $t\to \infty$ implying
\BQNY
 h(N(t)) &\toprob & \pk{(\kal{X}_1, \kal{Y}_1)  \le  \vk{w}_1 \ldot   (\kal{X}_m, \kal{Y}_m ) \le \vk{w}_m}, \quad
 t\to \infty.
\EQNY
If further $\vk{u}$ is such that $\pk{ \vk{U}_{\kal{Z}}= \vk{u}}=0$ then
$$\vk{1}( \vk{U}_{t,N(t)} \le \vk{u}) \todis   \vk{1}( \vk{U}_{\kal{Z}} \le \vk{u}), \quad t\to \infty, $$
consequently by Slutsky lemma (see e.g.\ Kallenberg (1997))
\BQNY
\vk{1}( \vk{U}_{t,N(t)}) h(N(t))& \todis & \vk{1}( \vk{U}_{\kal{Z}} \le \vk{u})
\pk{(\kal{X}_1, \kal{Y}_1)  \le  \vk{w}_1 \ldot   (\kal{X}_m, \kal{Y}_m ) \le \vk{w}_m}, \quad t\to \infty.
\EQNY
}
Since \ccr{$\vk{1}( \vk{U}_{t,N(t)}) h(N(t)),t>0$ is positive} and bounded (consequently uniformly integrable)
we have taking the expectation with respect to $N(t)$ and passing to limit
$$
\Bigl(   \vk{U}_{t, N(t)}, \vk{T}_{1}(N(t))\ldot \vk{T}_{m}(N(t)) \Bigr)
 \todis
\Bigl( \vk{U}_{\kal{Z}}, (\kal{X}_1, \kal{Y}_1) \ldot (\kal{X}_m, \kal{Y}_m ) \Bigr), \quad t\to \infty.
$$
Now, for any $ j\le m, t>0$ we may write
\BQNY
\lefteqn{\Biggl( \biggl(\frac{X_{N(t)-j+1:N(t)}- \bay  }{\aay} ,
\frac{Y_{N(t)-j+1:N(t)}- \bby }{\aby}\biggr) \Biggr)}\\
& =&
\Bigl( \tilde{a}_1(t,Z_t) T_{j1}(N(t))+ \tilde{b}_1(t,Z_t) ,
\tilde{a}_2(t, Z_t)T_{j2}(N(t))+ \tilde{b}_2(t,Z_t) \Bigr),
\EQNY
hence the proof follows by the continuous mapping theorem.
 \QED

\COM{
We may write for any $t$ positive and any fixed $j \le m$
\BQNY
\lefteqn{\Biggl( \biggl(\frac{X_{N(t)-j+1:N(t)}- \bay  }{\aay} ,
\frac{Y_{N(t)-j+1:N(t)}- \bby }{\aby}\biggr) \Biggr)}
\notag\\
& =&
\Biggl(\tilde{a}_1(t)\fracl{X_{N(t)-j+1:N(t)}- \cgb{b_1(N(t))} }{a_1(N(t)))} , \
\tilde{a}_2(t) \fracl{Y_{N(t)-j+1:N(t)}- b_2(N(t)) }{a_2(N(t))} \Biggr)\notag\\
&&+ \Biggl(\biggl(\frac{b_1([N(t)/t] t)- \bay}{\aay} ,
\frac{b_2([N(t)/t]t)-\bby }{\aby}\biggr) \Biggr),
\EQNY
where
$$
\tilde{a}_i(t):= \frac{a_i( [N(t)/t]t )}{a_i(t)}, \quad i=1,2.$$
Hence in view of \eqref{eq:1} and the assumptions on the counting process $N(t)$,
the joint weak convergence follows easily using further
Lemma 2.2.6 of Mikosch (2004) and Theorem 4.5 of Billingsley (1968). \QED
}
{\bf Proof of Proposition \ref{prop:2}:}
By the assumptions and \neprop{prop:1} we have
that \eqref{proof1} holds, hence  we may write for $t>0, p,q\inn$ using the
continuous mapping theorem 
\BQNY \label{eq:RR}
\lefteqn{
\Biggl( \frac{\stpX-  \bay c_1 }{\aay},
\frac{\stpY-  \bby c_2  }{\aby}\Biggr)
}\notag \\
&=&
\Biggl( \sum_{j=1}^{p}k_{j1}[X_{N(t)-j+1:N(t)}- \bay]/\aay,
 \sum_{j=1}^{q}k_{j2}[Y_{N(t)-j+1:N(t)}- \bby]/\aby \Biggr)\HR
&\todis & \Biggl( \kZ^{\gamma_{1}}\sum_{j=1}^{p} k_{j1}  \kal{X}_j+c_1\delta_{1} ,
\kZ^{\gamma_{2}}\sum_{j=1}^{q}k_{j2}\kal{Y}_j +c_2\delta_{2}\Biggr),
\quad t \to \IF,
\EQNY
with $\delta_{i}, \gamma_{i}, i=1,2$ as in \neprop{prop:1}, thus the proof is complete. \QED

{\bf Proof of Proposition \ref{prop:3}:} The proof follows immediately using  \neprop{prop:2} and
the result of Example 3 recalling further the expression for the joint density function given in \eqref{eq:hm}.
 \QED

\bigskip
{\bf Acknowledgment:} I would like to thank two  Referees for
several comments and in particular for correcting the result of
Example 3 pointing out a missing term related to Euler's constant.

\end{document}